\documentclass[10pt]{article}
\usepackage[english]{babel}
\usepackage{latexsym}
\usepackage{color}
\usepackage{amsmath}
\usepackage[latin1]{inputenc}
\usepackage{amsfonts}
\usepackage{amssymb}
\usepackage{mathrsfs}
\usepackage{eucal}
\usepackage{exscale}
\usepackage{makeidx}
\usepackage{makeidx}
\usepackage{graphicx}
\usepackage[T1]{fontenc}
\usepackage{amsmath,amssymb}
\usepackage[all]{xy}

\newtheorem{theorem}{Theorem}
\newtheorem{proposition}{Proposition}
\newtheorem{remark}{Remark}
\newtheorem{definition}{Definition}

\newtheorem{corollary}{Corollary}
\newtheorem{lemma}{Lemma}

\newtheorem{example}{Example}

\newcommand{\preuve}{\indent {\it Proof.}\hspace{4mm}}

\title{Almost positively closed models}
\author{Mohammed Belkasmi}
\begin{document}
\maketitle

\begin{abstract}
We introduce the notions of  almost positively closed 
models and positive strong amalgamation property. We study
the fundamental properties of  these notions
and develop  some interactions 
between them.
\end{abstract}
\section*{Introduction}
 The positive model theory is concerned essentially with the study
  of h-inductive theories which are built without the use of the
   negation.  Considering
 positive formulas instead of  formulas, and  homomorphisms instead of embeddings, the positive logic creates a new situations outside
  the framework of logic with negation. 
 In this paper, we will explore one of the specific aspects of positive logic which embodies the notions of  
 of algebraic 
closedness and strong amalgamation, and undertake to study  some interactions 
between these new notions inspired directly or indirectly from the 
works P.Bacsish \cite{bachich}.
 In the first section, after  summarising  the necessary background of positive
model theory, we introduce the general form of symmetric and 
asymmetric amalgamations. We show that the model-completeness
of an h-inductive theory can be characterized by a form of
 symmetric amalgamation. The second section is devoted to the 
 notions of almost positively closed models and a special class
 of positive formulas called $(A, T)$-closed formulas. Note that 
 the terminology "closed formula" here has different meaning of 
 the notion of formulas without free variables.
 We analyse the class of almost positively closed model
 and present a characterization through  some properties of the class of the $(A, T)$-closed formulas. In the third section, we introduce the notions of positive 
 strong amalgamation and h-strong amalgamation properties.
 We show that the class of almost positively closed has 
 the positive strong amalgamation property. Further
  we give a syntactic characterization of positive 
 strong amalgamation bases.

 \section{Positive model theory}
 \subsection{Basic definitions and notations}
 In this subsection we briefly introduce the basic terminology related 
to the positive logic. For more details, The 
reader is referred to \cite{ana}, \cite{poizat}. 
 
 Let $L$  be a first order language that contains the  symbol of equality and a constant $\perp$ denoting the antilogy. 
 The quantifier-free positive formulas are obtained from atomic formulas using the connectives $\wedge$ and $\vee$. The positive formulas are build from quantifier-free positive formulas using 
the logical operators and quantifier  $\wedge, \vee$ and 
$\exists$ respectively. Eventually, the positive formulas are of 
the form $\exists\bar y\ \phi(\bar x, \bar y)$, where 
$\phi(\bar x, \bar y)$ is quantifier-free formula. The variables 
$\bar x$ are said to be free. And a sentence is a formula without 
free variables.\\
  A sentence is said to be h-inductive, if it is a finite conjunction of sentences of the form: 
  $$\forall\bar x\ (\varphi(\bar x)\rightarrow\psi(\bar x))$$
   where $\varphi(\bar x)$ and $\psi(\bar x)$ are positive formulas. The h-universal sentences are the negation of 
   positive sentences.
 
 Let $A$ and $B$ be two $L$-structures. A map $h$ 
from $A$ to $B$  is  a homomorphism if  
for every  tuple $\bar a\in A$ 
and for every atomic formula $\phi$, if $A\models\phi(\bar a)$
then $B\models\phi(h(\bar a))$.
So we say that $B$ is  a continuation of $A$.\\
An embedding of $A$ into $B$ is a homomorphism 
$h:\ A\longrightarrow B$ such that,  for every $\bar a\in A$
and $\phi$ an 
atomic formula, if $B\models\phi(h(\bar a))$
then  $A\models\phi(\bar a)$. The homomorphism 
$h:\ A\longrightarrow B$ is said to be 
  an immersion whenever for every  $\bar a\in A$, 
 and $\phi$ a 
positive formula , if $B\models\phi(h(\bar a))$
then  $A\models\phi(\bar a)$.

A class of $L$-structures is said to be h-inductive if it is 
closed under the inductive limit of homomorphisms.

parallel to the role of existentially closed  structures in the framework of logic with negation, for every h-inductive theory
$T$ there exists a class of models of $T$ which represent the theory marvellously, and which enjoy the properties desired by every structures; namely, the h-inductive property of the class, the maximality of types (positive formulas satisfied by an element), amalgamation property, and others. These modules are called 
positively closed.
\begin{definition}
A model $A$ of an h-inductive theory $T$ is said to be 
positively closed (in short, pc)
if every homomorphism from $A$ into $B$ a model of $T$ is 
an immersion.
\end{definition}
The following lemmas announce the main properties of pc models. They will be used without mention.
\begin{lemma}[Lemma14, \cite{poizat}]
A model $A$ of an h-inductive theory is pc if and only if 
for every positive formula $\varphi$ and $\bar a\in A$, if 
$A\not\models\varphi(\bar a)$ then there exists a positive formula 
$\psi$ such that 
$T\vdash\neg\exists\bar x (\varphi(\bar x)\wedge\psi(\bar x))$ and $A\models\psi(\bar a)$.
\end{lemma}
\begin{lemma}[theorem 2, \cite{poizat}]
Every model on an h-inductive theory $T$ is continued in a pc 
model of $T$.
\end{lemma}
For every positive formula $\varphi$ we denote by 
$Ctr_T(\varphi)$ the set of positive formulas $\psi$ such that: 
$$T\vdash\neg\exists\bar x (\varphi(\bar x)\wedge\psi(\bar x)).$$
Two h-inductive theories are said to be companion  
if every model of one of them can be continued into a model of
the other, or equivalently if  the theories have the same pc models.\\ 
 Every $h$-inductive 
theory $T$ has a maximal companion denoted $T_k(T)$,  called
the Kaiser's hull of $T$.  $T_k(T)$
is the set of $h$-inductive sentences satisfied 
in each   pc models of $T$.
 Likewise,
$T$ has a minimal companion denoted $T_u(T)$, formed by its
$h$-universal consequence sentences. 
\begin{remark}
 Let $T_1$ and $T_2$  two h-inductive theories. 
The following  propositions are equivalent:
\begin{itemize}
\item $T_1$ and $T_2$ are companion.
\item $T_k(T_1)\equiv T_k(T_2)$.
\item $T_u(T_1)\equiv T_u(T_2)$.
\end{itemize}
\end{remark}
\begin{definition}
Let $T$ be an h-inductive theory.
\begin{itemize}
 \item $T$ is said to be model-complete 
if  every model of $T$ is a pc model of $T$.
\item We say that $T$ has a model-companion whenever $T_k(T)$ is 
model-complete theory.
\end{itemize}
 \end{definition}

 Let $A$ be a model of $T$. We shall use the following notations:
\begin{itemize}
\item $Diag^+(A)$,  the set  of  positive quantifier-free sentences satisfied by $A$ in the language $L(A)$. 
\item $Diag(A)$,  the set  of atomic and negated atomic  sentences satisfied by $A$ in the language $L(A)$. 
\item We denote by 
$T^+(A)$ the $L(A)$-theory $T\cup\{Diag^+(A)\}$.
\item $T_i(A)$ (resp. $T_u(A)$) denote the set of h-inductive
 (resp. h-universal ) 
 $L(A)$-sentences satisfied by $A$.
\item $T_i^\star(A)$ (resp. $T_u^\star(A)$) 
denote the set of h-inductive (resp. h-universal) $L$-sentences 
 satisfied by $A$.
 \item $T_k(A)$ (resp. $T_k^\star(A)$) denote the Kaiser's
hull of $T_i(A)$
(resp. of  $T_i^\star(A)$).
\item  For every subset $B$ of $A$, we denote  by $T_i(A, B)$ (resp. $T_u(A, B)$) the set of h-inductive (resp. h-universal ) $L(B)$-sentences satisfied by A.
\end{itemize}

\begin{definition}
Let $A$ and $B$ be two $L$-structures and $h$ a homomorphism  
from $A$ into $B$. 
 $h$ is said to be a strong immersion 
(in short s-immersion) if $h$ is an immersion and 
$B$ is a model of $T_i(A)$.
\end{definition}
\begin{remark}\label{remark3}
Let $A$ and $B$ two $L$-structures. 
We have the following properties:
\begin{enumerate}
\item If $A$ is immersed in $B$ then $T_u^\star(A)= T_u^\star(B)$,
and $T_i^\star(B)\subseteq T_i^\star(A)$.
\item $A$ is immersed in $B$ if and only if $T_i(B, A)\subseteq T_i(A)$.
\item $A$ is strongly immersed in $B$ if and only if 
$T_i(B, A)= T_i(A)$.
\item If $A$ and $B$ are two pc models of $T$ then every homomorphism
from $A$ into $B$ is a strong immersion. Indeed, let 
$\varphi(\bar a, \bar x)$  and $\psi(\bar a, \bar x)$ be two 
positive formulas and let $\chi$ the h-inductive sentence
$\forall\bar x (\varphi(\bar a, \bar x)\rightarrow
\psi(\bar a, \bar x))$. Suppose that 
$A\vdash\chi$ and $B\not\vdash\chi$, then there is $\bar b\in B$ such
that $B\models\varphi(\bar a, \bar b)$ and 
$B\not\models\psi(\bar a, \bar b)$. Given that $B$ is a pc model,
there exists $\psi'(\bar x, \bar y)\in Ctr_T(\psi(\bar x, \bar y))$
such that $B\models\psi'(\bar a, \bar b)$. Since $A$ is immersed 
in $B$, then there is $\bar a'\in A $ such that 
$A\models\varphi(\bar a, \bar a')\wedge\psi'(\bar a, \bar a')$,
which implies  $A\models\varphi(\bar a, \bar a')$ and 
$A\not\models\psi(\bar a, \bar a')$, contradiction.
\item The pc models of the $L(A)$-theory $T^+(A)$   are
the pc models of $T$ that are continuation of $A$. Indeed, 
it is clear that every pc model of $T$ in which $A$ is continued 
is a model of $T^+(A)$  and then a pc model of $T^+(A)$.\\
Conversely, let 
$B$ be a pc model $T^+(A)$ and $C$ a pc model of $T$ in which 
$B$ is continued by a homomorphism $f$.  Then  $C$ is a continuation of $A$, so  $C$ is a model of $T^+(A)$. Thereby $f$ is an immersion, which implies that $B$ is a pc model of $T$.
\end{enumerate}
\end{remark}
Let $A$ and $B$ two be $L$-structures and 
$f$ a mapping from $A$ into $B$. 
We will use the following notations: 
\begin{itemize}
\item $Hom(A, B)$ the set homomorphisms from
$A$ into $B$.
\item $Emb(A, B)$ the set embeddings from
$A$ into $B$.
\item $Imm(A, B)$ the set immersions from
$A$ into $B$.
\item $Sm(A, B)$ the set s-immersions from
$A$ into $B$.
\end{itemize} 
\begin{remark}
Let $A$ and $B$ be two $L$-structures and 
$f$ a mapping from $A$ into $B$. Consider $B$ 
as a $L(A)$-structure by interpreting 
the elements of $A$ in $B$ by  $f$. We have the 
following: 
\begin{itemize}
\item $f\in Hom(A, B)$ if and only if $B$ is a model of 
$Diag^+(A)$. 
\item $f\in Emb(A, B)$ if and only if $B$ is a model of 
$Diag(A)$.
\item $f\in Imm(A, B)$ if and only if $B$ is a model of 
$Diag^+(A)\cup T_u(A)$.
\item $f\in Sm(A, B)$ if and only if $B$ is a model of 
$T_i(A)$.
\end{itemize}
\end{remark}
\subsection{Positive amalgamation}
To abbreviate the nominations of homomorphism , embedding, 
immersion and strong immersion in the definition
of the notions of amalgamation,
we will use the first letter
of each mapping defined above.
\begin{definition}\label{dfamalgamation}
 Let $\Gamma$ be a class of $L$-structures and $A$ a member
  of $\Gamma$.
We say that $A$ is 
$[h, e, i, s]$-amalgamation basis 
of $\Gamma$, if for every $B, C$ members of $\Gamma$, if 
$A$ is continued into $B$ by $f$ and embedded into $C$ by 
$g$,  there  exist
$D\in \Gamma$, $g'\in Imm(C, D)$  and $f'\in Sm(B, D)$
 such that 
the following diagram commutes:
\[
\xymatrix{
    A \ar[r]^{g} \ar[d]_{f} & {C} \ar[d]^{g'} \\
    B \ar[r]_{f'} & {D}
  }
\]
By the same way we define the notion of 
$[\alpha, \beta, \gamma, \delta]$-amalgamation property for every 
 $(\alpha, \beta, \gamma, \delta)\in \{h, e, i, s\}^4$.\\
We say that $A$ is an 
$[\alpha]$-amalgamation basis of $\Gamma$,
if $A$ is an 
$[\alpha, \alpha, \alpha, \alpha]$-amalgamation basis 
of $\Gamma$.\\
We say that $A$ is $[\alpha, \beta]$-symmetric 
amalgamation basis  of $\Gamma$ whenever $A$ is an 
 $[\alpha, \beta, \beta, \alpha]$-amalgamation basis of 
 $\Gamma$.\\
We say that $A$ is $[\alpha, \gamma]$-asymmetric 
amalgamation basis  of $\Gamma$, whenever $A$ is an 
 $[\alpha, \beta, \alpha, \beta]$-amalgamation basis of 
 $\Gamma$.\\
\end{definition}
The following remark list some properties of diver forms of 
amalgamation  with the notations and terminology given in the 
definition above.
\begin{remark}\label{examplamalgam}
\begin{enumerate}
\item Every $L$-structure $A$ is an $[i,h,s,h]$-amalgamation 
basis in the class of $L$-structures
(lemma 4,
 \cite{ana}). Since every strong immersion is an 
 immersion, it follows that every 
 $L$-structure $A$ is an $[s,h]$-asymmetric 
 amalgamation basis in the class of $L$-structures.
 \item Every $L$-structure $A$ is an $[s,i]$-asymmetric
  amalgamation basis in the class of $L$-structures
(lemma 5,
 \cite{ana}).
 \item Every $L$-structure $A$ is an $[e,s]$-asymmetric
  amalgamation basis in the class of $L$-structures
(lemma 4,
 \cite{ana2}).
\item Every $L$-structure $A$ is an $[i,h]$-asymmetric
 amalgamation basis in the class of $L$-structures
(lemma 8,
 \cite{poizat}).
 \item Every pc model of an h-inductive theory $T$ is an 
 $[h]$-amalgamation basis in the class of model of $T$.
\end{enumerate}
\end{remark}
\begin{lemma}
Every $L$-structure is $[s, x]$-asymmetric
amalgamation basis in the class of $L$-structure, where $x$
is a homomorphism, an embedding or an immersion.
\end{lemma}
\preuve 
All cases are given in  Remark \ref{examplamalgam}.
\begin{lemma} \label{hi sym}
A model of $T$ is  pc if and only if it has the  $[h, i]$- symmetric amalgamation property in the class of models of $T$.
\end{lemma}
\preuve 
Let $A$ be an $[h, i]$- symmetric amalgamation basis of $T$. 
Assume that  $A\not\models\varphi(\bar a)$, where $\bar a\in A$ and 
$\varphi$ a positive formula. 
Given that $A$ is an $[h, i]$- symmetric amalgamation basis,
we claim that $T\cup Diag^+(A)\cup\{\varphi(\bar a)\}$  is 
inconsistent. Otherwise, we can find two continuations of $A$
one of them satisfies $\varphi(\bar a)$ 
and the other does not satisfy $\varphi(\bar a)$, 
which contradicts the assumption that $A$ has the 
 $[h, i]$- symmetric amalgamation property.
 \begin{proposition}
 An h-inductive theory $T$ has a model companion
 if and only if $T_k(T)$ has the $[h, i]$- symmetric amalgamation
 property.
 \end{proposition}
 \preuve 
 Suppose that $T$ has a model companion, then every model 
 of $T_k(T)$ is a pc model. Since the pc models have the 
 $[h]$-amalgamation property and the homomorphisms between 
 the pc models are immersion, it follows from the fifth bullet 
 of the remark \ref{examplamalgam} that 
$T_k(T)$ has the $[h, i]$- symmetric amalgamation
 property.\\
 The opposite direction follows easily from the lemma \ref{hi sym}.
 
\section{Almost positively closed structures}
In this section, we introduce the notions of  almost 
and $\Gamma$-almost positively closed models, we give 
a syntactic characterisation and a characterization
via  the  closed formulas  
which turns out to be an essential tool in the the study of the notion 
of $\Delta$-almost positively closedness.
\begin{definition}
Let $T$ be an h-inductive theory, and let 
 $\Delta$ be a subset of $L$-quantifier-free 
 positive formulas.
  A  model $A$ of $T$ is said to be: 
  \begin{itemize}
  \item Almost positively closed (apc in short),
 if for every model 
$B\models T$,  $f\in Hom(A, B)$ and 
$\varphi(\bar x, \bar y)$  a quantifier-free positive formula,
 if $B\models\exists\bar y\varphi(\bar a, \bar y)$
  and $\bar a\in A$
 then there is 
$\bar a'\in A$ such that $B\models\varphi(\bar a, \bar a')$.
\item $\Delta$-almost positively closed ($\Delta$-apc in short),
 if for every model 
$B\models T$,  $f\in Hom(A, B)$ and 
 $\varphi(\bar x, \bar y)\in\Gamma$, 
 if $B\models\exists\bar y\varphi(\bar a, \bar y)$
  and $\bar a\in A$
 then there is 
$\bar a'\in A$ such that $B\models\varphi(\bar a, \bar a')$.
\item Weakly almost positively closed (wpc in short),
 if for every pc model 
$B\models T$,  $f\in Hom(A, B)$ and 
$\varphi(\bar x, \bar y)$  a quantifier-free positive formula
 if $B\models\exists\bar y\varphi(\bar a, \bar y)$
  and $\bar a\in A$
 then there is 
$\bar a'\in A$ such that $B\models\varphi(\bar a, \bar a')$.
\item 
 $\Delta$-weakly almost positively closed 
($\Delta$-wpc in short),
 if for every pc model 
$B\models T$,  $f\in Hom(A, B)$ and 
 $\varphi(\bar x, \bar y)\in\Gamma$),
 if $B\models\exists\bar x\varphi(\bar a, \bar x)$
 then there is 
$\bar a'\in A$ such that $B\models\varphi(\bar a, \bar a')$.
\end{itemize}   
\end{definition}
 \begin{theorem}\label{charactirization apc model}
 Let $A$ be a model of an h-inductive $L$-theory $T$,  
 let $\Delta$ be a set of $L(A)$-quantifier free positive 
 formula. $A$ is $\Delta$-apc of $T$ if and only if for every 
 $\varphi(\bar a, \bar x)\in \Gamma$, there exists a 
 quantifier free positive formula $\psi(\bar a, \bar a')\in 
 Diag^+(A)$ such 
 that 
 $$T\vdash\forall\bar x\bar{y}((\psi(\bar x,\bar y)\wedge\exists\bar z
 \varphi(\bar x, \bar z))\rightarrow\varphi(\bar x,\bar y)).$$ 
 \end{theorem}
 \preuve 
 Assume that $A$ is an $\Delta$-apc model of $T$ and let 
 $\varphi(\bar a, \bar x)\in \Gamma$. Considering  
 $T^\ast= 
 T\cup Diag^+(A)\cup \{\exists\bar x\varphi(\bar a, \bar x)\}$
 is consistent and $A$ is $\Delta$-apc then 
 $T^\ast\cup\{\neg\varphi(\bar a, \bar a')|\bar a'\in A\}$
 is inconsistent. Thereby there are $\bar a'\in A$ and 
 $\psi(\bar a, \bar a')\in Diag^+(A)$ such that 
 $T\cup\{\psi(\bar a, \bar a'), \neg\varphi(\bar a, \bar a'), 
 \exists\bar x\varphi(\bar a, \bar x)\}$ is inconsistent, which 
 implies 
 $T\vdash \forall\bar x, \bar y((\psi(\bar x, \bar y)\wedge
 \exists\bar z\varphi(\bar x, \bar z))\rightarrow\varphi(\bar x,
 \bar{y}))$.\\
 The other direction is clear.
 \begin{corollary}\label{lemma theorem 1}
 Let $A$ be a model of $T$ and $\Delta$ a set of quantifier-free 
 positive $L(A)$-formulas. If $A$ is immersed in an $\Delta'$-apc
 model $B$ of $T$ and $\Delta\subseteq \Delta'$, then 
 $A$ is an $\Delta$-apc model of $T$.
 \end{corollary}
 \preuve 
 Let $\varphi(\bar a, \bar x)\in\Delta$. Given that 
 $\varphi(\bar a, \bar x)\in \Delta'$ and $B$ is $\Delta'$-apc
  model of $T$, by  theorem \ref{charactirization apc model}
there exists 
$\psi(\bar a, \bar b)\in Diag^+(B)$ such that 
$$T\vdash \forall\bar x, \bar y((\psi(\bar x, \bar y)\wedge
 \exists\bar z\varphi(\bar x, \bar z))\rightarrow\varphi(\bar x,
 \bar{y})).$$
 On the other hand, since $A$ is immersed in $B$ then there is 
 $\bar a'\in A$ such that $\psi(\bar a, \bar a')\in Diag^+(B)$, 
 hence  $A$ is an $\Delta$-apc of $T$ by theorem
  \ref{charactirization apc model}.
\begin{remark}\label{remarksapcmodel}
Let $T$ be an h-inductive $L$-theory and $\Delta$ a set of 
quantifier-free positive $L$-formulas. We have the following 
properties:
\begin{enumerate}
\item If $A$ is apc then $A$ is wpc of $T$. 
\item Every pc model of $T$ is an apc (resp. $\Delta$-apc) 
model of $T$.
\item The classes of apc and wpc (resp. $\Delta$-apc and 
$\Delta$-wpc) models of $T$ are $h$-inductive.
\item If $A$ is an apc model of $T$ and 
  $B$ a model of $T$, then $Emb(A, B)=Imm(A, B)$.
  \item Let $\Delta\subseteq\Delta$ be two sets of free 
  quantifier positive formulas. If $A$ is $\Delta$-apc 
  (resp.  $\Delta$-wpc) then $A$ is $\Delta$-apc 
  (resp.  $\Delta$-wpc).
\item A model $A$ of $T$ is  apc if and only if for every 
positive quantifier-free formulas $\varphi(\bar x, \bar y)$, there 
exists $\psi(\bar a, \bar a')\in Diag^+(A)$ such that 
$$T\vdash\forall\bar x\bar y((\psi(\bar x, \bar y )\wedge\exists\bar z\varphi(\bar x, \bar z))\rightarrow\varphi(\bar x, \bar y)).$$
\item Every apc model of $T$ has the property of 
  $[e, h]$-asymmetric amalgamation( property $4$ of Remark
 \ref{remarksapcmodel}, and the property $4$ of Remark 
 \ref{examplamalgam}). 
\end{enumerate}
\end{remark}
\begin{example}
\begin{enumerate}
\item Let $L=\{f\}$ functional language. Let $T$ be the h-inductive 
theory $\forall x, y (f(x)=f(y)\rightarrow x=y)$. The theory $T$ has a model companion axiomatized by $T_k(T)=T\cup \{\forall xy\ 
(x=y)\}$. The class of apc model of $T$ is elementary and axiomatized by the h-inductive theory
$$T\cup\{\exists x, f(x)=x\}\cup\{\forall x\exists y (f(y)=x)\}.$$
\item Let $L$  and $T$ respectively the functional language  and 
the theory defined above. Let $T''$ the h-inductive theory 
$T\cup \{\neg\exists x\ (f(x)=x)\}$. The class of apc model of 
$T''$ is axiomatized by the h-inductive theory 
$$T\cup\{\forall x\exists y (f(y)=x)\}\cup\{\exists x f^p(x)=x|\ 
p\ \text{prime number}\}.$$
\item Let $T_f$ the theory of fields. Since the negation of equality $x=y$ is defined by the positive formula 
$\exists z\ (x-y)\cdot z=1 $ and every homomorphism is an embedding 
then the classes of apc fields, pc fields and existentially closed 
fields are equals.
\end{enumerate}
\end{example}
\begin{definition}\label{closedformula}
Let $T$ be an h-inductive $L$-theory  and $A$ a model of $T$.
\begin{itemize}
\item A positive formula $\varphi(x)$ is said to be 
$T$-algebraic if there exists a positive formula 
$\psi(y_1,\cdots, y_n)$ such that $\psi\not\equiv \perp$ modulo
$T$ (ie, $\psi(y_1,\cdots, y_n)$ has a realisation in some
model of $T$) and;
$$
T\vdash\forall x,\bar y
((\varphi( x)\wedge
 \psi(\bar y))\rightarrow 
 \bigvee_{i}x=y_i).$$
 We denote by $Al_T$ the set of 
 $T$-algebraic quantifier-free positive $L$-formulas.
 \item  For every 
 formula $\varphi(\bar x)$, we denote by
  $E(\varphi, T)$
 the set of positive formulas $\psi(\bar y)$
 such that $\psi\not\equiv \perp$ modulo
$T$
  that satisfy 
 the property:
 $$T\vdash\forall\bar x\bar y
((\varphi(\bar x)\wedge
 \psi(\bar y))\rightarrow 
 \bigvee_{i, j}x_i=y_j))$$
\item 
A positive formula $\varphi(\bar x, \bar y)$ is said to be
$(A, T)$-closed if 
for every pc model continuation $B$ of $A$, if 
$B\models\varphi(\bar a, \bar b)$ for some $\bar a\in A$
then $\bar b\in A$.
\end{itemize}
\end{definition}
\begin{remark}\label{remark6} 
\begin{enumerate}
\item A quantifier-free positive formula is $T$-algebraic
if and only if its is algebraic  in the sense of Robinson \cite{robinson}.
 \item 
 Given that the class of pc models of $T^+(A)$ coincides 
 with the class of pc models of $T$ that are continuation of 
 $A$ (bullet 5 of  remark \ref{remark3}), then 
 a formula is 
$(A, T)$-closed if and only if it is $(A, T^+(A))$-closed.
\item Let $A$ be a model of $T$. Denote by $C_A$ the set 
of quantifier-free formulas that are  $(A, T)$-closed. By definition of $C_A$ we observe that 
$A$ is $C_A$-wpc.
\item  If 
every formula in $Al_{T^+(A)}$ is $(A, T)$-closed, by the property 
2  above and the definition of formulas $(A, T)$-closed,  we conclude that   
$A$ is $Al_{T^+(A)}$-wpc. Considering that 
$Al_T\subset Al_{T^+(A)}$, $A$ is also $Al_T$-wpc.
\end{enumerate}
\end{remark}
\begin{definition}
For every quantifier-free
positive formulas $\varphi(\bar x)$ such that $\varphi(\bar x)\not\equiv\perp$ modulo $T$, we denote by $E_T$ the set of 
 quantifier-free
positive formula $\varphi(\bar x)$ such that 
$E(\varphi, T)\neq\emptyset$.
\end{definition}
\begin{lemma}\label{carac amalg apc}
  Let $A$ be an  h-amalgamation basis  of $T$.
  If $A$ is $E_{T^+(A)}$-wpc (resp. $E_{T^+(A)}$-apc) then
every  formula in $Al_{T^+(A)}$ is $(A, T)$-closed.    
\end{lemma}
\preuve  Let $A$ be 
an $E_{T^+(A)}$-wpc and an h-amalgamation basis of $T$.
Assume the existence of a formula  
 $\varphi(\bar a,  y)\in Al_{T^+(A)}$ 
 such that $\varphi(\bar a,  y)$ is not $(A, T)$-closed.
 So, there exist $B$ a pc models of $T$ and $b\in B-A$ 
 such that  $B\models\varphi(\bar a,  b)$. Let  
 $\psi(\bar a, \bar x)\in E(\varphi, T^+(A))$. Let $C$ be  a
 pc model of $T^+(A)$ and $\bar c\in C$ such that 
 $C\models\psi(\bar a, \bar c)$. 
  Given that $\psi(\bar a, \bar y)\in E_{T^+(A)}$ and $A$ is an 
  $E_{T^+(A)}$-wpc model of $T$, then  there is $\bar a'$ in 
  $A$ such that $C\models\psi(\bar a, \bar a')$. 
 Let  $D$ be a model of $T$ that 
   amalgamate commutatively the diagram 
   $C\leftarrow A\rightarrow B$, so $B$ is immersed 
   in $D$ and  
   $B\models \varphi(\bar a,  b)\wedge\psi(\bar a, \bar a')$.
  Which implies $\bigvee_{i}b= a_i$, contradiction.\\
The proof of the case where $A$ is  $Al_{T}$-apc is an immediate  
consequence of this Lemma.
  
\section{Strong amalgamation}
In this section we introduce 
 the notions of positive strong amalgamation and 
h-strong amalgamation. We investigate their properties and 
 interactions with the notion of  almost positively closedness.
\subsection{positive strong amalgamation}
 \begin{definition}\label{dfamalgamation} 
 Let $T$ be an h-inductive theory. A  model $A$
  of $T$ is said to be 
  a positive strong amalgamation basis (in short PSA) 
  (resp.  h-strong amalgamation basis (in short h-SA))
of $T$, if for every pc models  (resp. models) $B$ and   $C$
of $T$, if 
$A$ is continued into $B$ and $C$ by two homomorphisms $f$ 
and $g$  respectively, then there  exist
$D$ a model of $T$, and   $f', g'$  two homomorphisms 
 such that 
the following diagram commutes:
\[
\xymatrix{
    A \ar[r]^{f} \ar[d]_{g} & {B} \ar[d]^{f'} \\
    C \ar[r]_{g'} & {D}
  }
\]
and satisfies the following property:\\
$\forall(b, c)\in B\times C$, if $g'(c)=f'(b)$ then 
there is $a\in A$ such that $c=g(a)$ and $b= f(a)$.\\
\end{definition}
\begin{example}
\begin{enumerate}
\item Let $T$ be an h-inductive theory such that a model $A$ of $T$ is pc 
if and only if $A\models\forall x,y\, x=y$. Then $T$ has the 
positive strong amalgamation property. As examples of these theories we have  the theory of groups and the theory of partially  ordered sets.
\end{enumerate}
\end{example}
\begin{lemma}\label{h-asymmetric stong}
Let $A, B, C$ be three $L$-structures. Let  $i\in Imm(A, B)$ and
$h\in Hom(A, C)$. Then there exist 
$D$ a $L$-structure,
 $h'$ a homomorphisms and $s$ an s-immersion,
 such that the following diagram commutes:
\[
\xymatrix{
    A \ar[r]^{i} \ar[d]_{h} & {B} \ar[d]^{h'} \\
    C \ar[r]_{s} & {D}
  }
\]
and satisfies the following property:\\
$\forall(b, c)\in B\times C$, if $h'(b)=s(c)$, then 
there exists $a\in A$ such that 
$c=h(a)$ and $b=i(a)$.
\end{lemma} 
\preuve   The proof
consists  in  the verification that 
 the set  
$$T_i(C)\cup Diag^+(B)\cup Diag^+(C)
\cup \{b\neq c|\ b\in B-A, c\in C-h(A)\}$$
is $L(B\cup C)$-consistent.\\
Assume, to the contrary that the set above is 
$L(B\cup C)$-inconsistent. Then there exist 
$\varphi(h(\bar a), \bar c)\in Diag^+(C)$,  
$\psi(\bar a, \bar b)\in Diag^+(B)$ where
$\bar{c}\in C-h(A)$  and $\bar{b}\in B-A$,  such that; 
$$T_i(C)\vdash \forall\bar y\ \ 
 ((\varphi(h(\bar a), \bar c)\wedge
 \psi(\bar a, \bar y))\rightarrow 
 \bigvee_{i, j}y_i=c_j).$$
Given that $B\models\psi(\bar a, \bar b)$ and 
$A$ is immersed in $B$, there is $\bar a'\in A$ such that 
$A\models\psi(\bar a, \bar a')$. Consequently,
$C\models\varphi(h(\bar a), \bar c)
\wedge\psi(h(\bar a),  h(\bar a')$.  Thereby
$C\models \bigvee_{i, j}h(\bar a')_i=c_j$,
which is a contradiction.
\begin{corollary}
Every   pc model $A$ of $T$ is a h-strong 
amalgamation basis of $T$.
\end{corollary}
\preuve 
Immediate from Lemma \ref{h-asymmetric stong}.
\begin{proposition}
Let $A$ and $B$ be two models of an h-inductive theory $T$, and 
let $i\in Imm(A, B)$. If $B$ is a h-SA basis of $T$ then $A$ is 
a PSA basis of $T$.
\end{proposition}
\preuve Let $A_1$ and $A_2$ be two pc mode of $T$,
$f\in Hom(A, A_1)$ and $g\in Hom(A, A_2)$. Considering 
the diagrams  $A_1\leftarrow A\rightarrow B$ and $A_1\leftarrow A\rightarrow B$,
 by  Lemma \ref{h-asymmetric stong}, we get the  commutative diagrams (1) and (2) below, where $f', g'$ are homomorphisms and 
 $i_1, i_2$ are strong immersions.\\
Now, given that  $B$ has the h-strong amalgamation property, we get 
the  commutative diagram (3) below:
  \[
\xymatrix{
A_1\ar[r]^{i_1}\ar@{}|{(1)}[rd]& B_1\ar[dr]^{f''}  &\\
 A \ar[r]^{i} \ar[d]_{g} \ar[u]^{f} &{B} \ar[d]^{g'}
  \ar[u]_{f'}\ar@{}|{(3)}[r]&C \\
A_2\ar[r]_{i_2}\ar@{}|{(2)}[ru]& B_2\ar[ur]_{g''} &    
  }
\]
where $ f'', g''$ are homomorphisms.\\
 We claim that $C$  makes the diagram
 $A_1\leftarrow A\rightarrow A_2$  strongly amalgamable. 
 Indeed, let $a_1\in A_1$ and 
$a_2\in A_2$ such that $f''\circ i_1(a_1)=g''\circ i_2(a_2)$, 
by the $h$-strong amalgamation property of the diagram (3), there
is $b\in B$ such that $f'(b)=i_1(a_1)$ and $g'(b)=i_2(a_2)$. 
Considering the properties of the diagrams $(1)$ and $(2)$, we 
get two  elements $a$ and $a'$ from $A$ such that:
$$
\left\{
    \begin{array}{ll}
        f(a)=a_1,& i(a)= b \\
        g(a')=a_2,& i(a')= b.
        
    \end{array}
\right.
$$ 
Given that $i$ is an immersion, then $a=a'$ and $f(a)=a_1,
  g(a)=a_2$. So $A$ is a PSA basis of $T$.

\begin{lemma}\label{lemma8}
Let $A$ be an h-amalgamation basis of $T$, 
$B$ a pc model of $T$ and     
$f\in Hom(A, B)$. $A$ is PSA basis of $T$ if and only 
if for every formula  $\varphi(\bar a, \bar x)\in E_{T^+(A)}$ 
and for every $b_1,\cdots, b_n\in B-f(A)$, we have 
$B\not\models\varphi(\bar a, b_1\cdots b_n)$.
\end{lemma}
\preuve 
Let $A$ be a PSA basis of $T$. Suppose that 
there are $\varphi(\bar a, \bar x)\in 
E_{T^+(A)}$, $B$ a pc model of $T$ and $f\in Hom(A, B)$ such that 
$B\models\varphi(\bar a, b_1\cdots b_n)$, where 
$b_1,\cdots, b_n\in B-f(A)$. Let 
$\psi(\bar a, \bar y)\in E(\varphi, T^+(A))$,  $C$ a 
 pc model of $T$ and $g\in Hom(A, C)$
 such that $C\models \psi(\bar a, \bar c)$.\\
 Given that $A$ is a PSA basis of $T$, we obtain the following
 commutative diagram 
\[
\xymatrix{
    A \ar[r]^{f} \ar[d]_{g} & {B} \ar[d]^{i} \\
    C \ar[r]_{i'} & {D}
  }
\]
where $i$ and $i'$ are immersion, and $D$  a model of $T$ 
 that satisfies the property: 
$$\forall(b, c)\in B\times C;\ i(b)=i'(c)\Rightarrow
\exists a\in A, (f(a)=b\wedge g(a)=c).$$
Now, since 
$D\models\varphi(\bar a, i(\bar b))\wedge\psi(\bar a, i'(\bar c))$
then $D\models \bigvee_{i, j}i(b_i)=i'(c_j)$, which implies the existence 
of an element $a'\in A$ such that $i(b_i)=i'(c_j)=i\circ f(a')$. Contradiction.

We shall prove the other direction by contrapositive.  Let 
$B$ and $C$ be two pc models of $T$,
$f\in Hom(A, B)$ and $g\in Hom(A, C)$
 such that the diagram 
\[
\xymatrix{
    C &A \ar[r]^{f} \ar[l]_{g} & B
  }
\]
is not h-strongly amalgamable. 
Thereby there exist 
$\varphi(f(\bar a),  b_1\cdots b_n)\in Diag^+(B)$,  
$\psi(g(\bar a), c_1\cdots c_m))\in Diag^+(C)$ where
$b_1,\cdots, b_n\in B-A$  and $c_1,\cdots, c_m\in C-A$,  such that 
$$T^+(A)\vdash \forall\bar y\ \ 
 ((\varphi(\bar a, \bar x)\wedge
 \psi(\bar a, \bar y))\rightarrow 
 \bigvee_{i, j}x_i=y_j).$$
Thus  $\varphi(\bar a, \bar x)\in E_{T^+(A)}$ 
and $B\models\varphi(\bar a,  b_1\cdots b_n)$.
\begin{theorem}\label{strongamalgam2}
 Let $A$ be an h-amalgamation basis of $T$, we have the following 
 properties:
 \begin{enumerate}
 \item If $A$ is a $E_{T^+(A)}$-wpc model of $T$  
  then $A$ is a PSA basis of $T$. 
  \item If $A$ is a PSA basis of $T$ then $A$ is $Al_{T^+(A)}$-wpc.
 \end{enumerate}
\end{theorem}
\preuve 
\begin{enumerate}
\item 
Let $A$ be an h-amalgamation 
basis and a $E_{T^+(A)}$-wpc model of $T$. Let $B$ and $C$ be two pc models of $T$, 
$f\in Hom(A, B)$ and $g\in Hom(A, C)$. Let 
 $D$ a model of $T$ such that the following 
 diagram commutes:
 \[
\xymatrix{
    A \ar[r]^{f} \ar[d]_{g} & {B} \ar[d]^{i_1} \\
    C \ar[r]_{i_2} & {D}
  }
\]
 where $i_1$ and $i_2$ are immersions.\\  
 We claim that  the   set 
 $T\cup Diag^+(B)\cup Diag^+(C)
\cup \{b\neq c|\ b\in B-A, c\in C-A\}$ is $L(B\cup C)$-consistent
(note that the element of $A$ are interpreted by the same symbols of constants in $B$ and $C$).
Suppose to the contrary  that the set above 
is inconsistent.  
Thus there are $\bar a\in A,\, \bar b\in B-A, \bar c\in C-A$,
$\varphi(\bar a, \bar b)\in Diag^+(B)$ and 
$\psi(\bar a, \bar c)\in Diag^+(C)$ such that 
 $$T\cup\{\varphi(\bar a, \bar b), 
 \psi(\bar a, \bar c), 
 \bigwedge_{i,j}b_i\neq c_j\}$$
 is $L(B\cup C)$-inconsistent, thereby 
 $$T^+(A)\vdash \forall\bar y,\bar z\ \ 
 ((\varphi(\bar a, \bar y)\wedge
 \psi(\bar a, \bar z))\rightarrow 
 \bigvee_{i, j}y_i=z_j).$$
 Now, since $C\models\psi(\bar a, \bar c)$,
  $\psi\in E_{T^+(A)}$, and 
 $A$ is an $E_{T^+(A)}$-wpc model, then there is 
 $\bar a'\in A$ such that 
 $C\models\psi(\bar a, \bar a')$. Thereby
 $D\models\psi(\bar a, \bar a')$, so  
 $B\models\psi(\bar a, \bar a')\wedge 
 \varphi(\bar a, \bar b)$.  
Which implies $B\models \bigvee_{i, j}b_i=a'_j$,
 contradiction. Thus
 $A$ is a PSA basis of $T$.
 \item   Suppose 
that $A$ is PSA of $T$. Since 
$ Al_{T^+(A)}\subseteq E_{T^+(A)}$, by Lemma \ref{lemma8}
every formula in $Al_{T^+(A)}$ is $(A,T)$-closed.  which implies 
that $A$ is a $Al_{T^+(A)}$-wpc model of $T$
by Remark \ref{remark6} (4).
\end{enumerate}
\bibliographystyle{plain}

\begin{thebibliography}{10}
\bibitem{ana}
Mohammed ~Belkasmi.
\newblock {\em Positive model theory and amalgamations}.
\newblock {\em Notre Dame Journal of Formal Logic}, vol. 55, (2014), 205-229.
\bibitem{ana2}
Mohammed ~Belkasmi.
\newblock {\em Positive amalgamation}.
\newblock {\em Logica Universalis}, Published online.
\bibitem{poizat}
Ita\"\i\ ~Ben Yaacov, Bruno ~Poizat.
\newblock {\em Fondements de la logique positive}.
\newblock {\em Journal of Symbolic Logic}, 72, 4, 1141--1162, 2007.
\bibitem{robinson}
Abraham ~Robinson.
\newblock {\em Introduction to model theory and the 
mathematics of algebras}.
\newblock {\em North Holland, Amsterdam}, 
1965.
\bibitem{bachich}
Paul D. ~Bacsich.
\newblock {\em Defining algebraic elements}.
\newblock {\em The Journal of Symbolic Logic}, 
Vol. 38, No. 1 (1973), pp. 93-101.
\end{thebibliography}

\begin{flushleft}
Department of Mathematics, College of sciences\\ 
 Qassim University\\
  Buraydah. Saudi Arabia
\end{flushleft}

\end{document}